\documentstyle[11pt]{article}

\newtheorem{theo}{Theorem}
\newtheorem{coro}{Corollary}
\newtheorem{prop}{Proposition}
\newtheorem{lemm}{Lemma}

\begin{document}

\def\ot{\otimes}
\def\we{\wedge}
\def\wec{\wedge\cdots\wedge}
\def\op{\oplus}
\def\ra{\rightarrow}
\def\lra{\longrightarrow}
\def\fso{\mathfrak so}
\def\cO{\mathcal{O}}
\def\cS{\mathcal{S}}
\def\fsl{\mathfrak{sl}}
\def\PP{\mathbb P}\def\PP{\mathbb P}\def\ZZ{\mathbb Z}\def\CC{\mathbb C}
\def\RR{\mathbb R}\def\HH{\mathbb H}\def\OO{\mathbb O}

\title{On the asymptotics of Kronecker coefficients}
\author{Laurent \sc{Manivel} \\  \\
UMI 3457 CNRS/Centre de Recherches Math\'ematiques,\\
Universit\'e de Montr\'eal, Canada\\
\texttt{manivel@math.cnrs.fr} }
\maketitle

\section{Introduction}

Kronecker coefficients are the structure constants for the tensor products of irreducible representations of symmetric groups. 
Their computation is thus an old and important problem in finite group theory. Through Schur-Weyl duality, they can be 
interpreted as multiplicities in Schur powers of tensor products and their computation becomes a problem in invariant theory. 
Applying the Borel-Weil theorem we can see  Schur powers of tensor products as spaces of sections of equivariant line bundles 
on flag varieties, and the problem can be interpreted as a question about Hamiltonian actions on symplectic manifolds. 
It can also be considered in relation with marginal problems in quantum information theory. Each of these perspectives 
gives access to interesting pieces of information about  Kronecker coefficients, which can be very hard to obtain, or 
even to guess, from a different perspective. But although we have these several points of views on the Kronecker 
coefficients, each revealing some part of their properties,  they remain quite mysterious, and very basic problems
seem completely out of reach at this moment. Let us mention a few important ones. 
\begin{enumerate} 
\item No combinatorial formula is known. We would like the Kronecker coefficients to count some combinatorial objects. 
Or points in polytopes. 
\item A few linear conditions for Kronecker coefficients to be non zero are known. We would like to know all of them. 
In the terminology of this paper, we would like to know the faces of the Kronecker polyhedra. 
\item These polyhedra being described, we would like to know how far we are from being able to decide whether a Kronecker coefficient is zero or not. This is a saturation type problem, related to Mulmuley's conjecture that this decision problem is in P. 
\item We would like to understand stretched Kronecker coefficients. They are given by certain quasipolynomials, which we would 
like to be able to compute efficiently.
\end{enumerate}
The analogue problems for Littlewood-Richardson coefficients have been solved. The Littlewood-Richardson rule gives a 
combinatorial recipe for their computation, which can be interpreted as a count of integral points in polytopes, once 
translated in terms of hives or honeycombs, for example. The full list of linear inequalities is known and gives a solution 
to the famous Horn problem (see eg \cite{fulton} for a survey). 
They are enough to decide whether a Littlewood-Richardson coefficients is zero or not, a decision 
problem which is definitely in P. Finally the stretched versions have been studied and are in fact given by polynomials, not just quasipolynomials \cite{rassart}. 
 
This gives some hope for the seemingly similar Kronecker coefficients. There are two different ways to appreciate the connections between these two families of numbers.  Littlewood-Richardson coefficients are multiplicities in Schur powers of direct sums rather than tensor products. They are also special Kronecker coefficients: in fact, quite remarkably there is a special facet of the Kronecker polyhedron which is exactly a Littlewood-Richardson polyhedron. The latter is known to be 
generated by the triple of partitions 
for which  the  Littlewood-Richardson coefficient is exactly one, and there is a remarquable equivalence
$$c_{\lambda,\mu}^\nu=1 \quad \Longleftrightarrow \quad c_{k\lambda,k\mu}^ {k\nu}=1 \; \forall k\ge 1.$$
Such a property will certainly not hold for Kronecker coefficients, but it suggests to introduce the following definition:

\medskip\noindent {\bf Definition.} 
A triple of partitions $(\alpha,\beta,\gamma)$ is {\it weakly stable} if the Kronecker coefficients 
$$ g(k\alpha,k\beta,k\gamma)=1 \quad  \forall k\ge 1.$$
It is {\it stable} if $ g(\alpha,\beta,\gamma)\ne 0$  and for any triple $(\lambda,\mu,\nu)$, the sequence of Kronecker coefficients 
$ g(\lambda+k\alpha,\mu+k\beta,\nu+k\gamma)$ is bounded (or equivalently, eventually constant). 

\medskip
Stability was introduced in \cite{stem} by Stembridge, who proved that it implies weak stability and conjectured that the two notions are in fact equivalent. From the Kronecker polyhedron perspective, stable or weakly stable  
triples will correspond to very special boundary points. 

One of the main objectives of the present paper is to explain how to construct large families of stable triples. Moreover we will
be able to describe the  Kronecker polyhedron locally around these special points, including by giving explicit equations of the facets they belong to. 

Our main tool for studying the Kronecker coefficients will be Taylor expansion. This might look strange at first sight, but 
recall that using Schur-Weyl duality and the Borel-Weil theorem we can interpret them in terms of  spaces of sections of
line bundles on flag manifolds. In this context it is very natural to use some basic analytic tools, like Taylor expansion, 
in order to analyze such sections. We will use these Taylor expansions around certain flag subvarieties, and the fact
that these subvarieties allow us to control the  Kronecker coefficients will depend on certain combinatorial properties
of standard tableaux. The key concept here is that of additivity, studied in detail by Vallejo \cite{vallejo, vallejo2}.
In fact this concept already appears in \cite{man-if}, where very similar ideas are introduced and used to understand
certain asymptotic properties of plethysms. Section 3 of \cite{man-if}  is in fact already devoted to Kronecker coefficients:
we explained how the method could be used in this context, what was the role of the additivity property, how we could deduce
infinite families of stable triples, and much more. In the present paper we explain our approach in greater detail in the 
specific context of Kronecker coefficients, and go a bit further than in \cite{man-if}. Our results are the following:
\begin{enumerate}     
\item We show that in a very general setting, stable triples can be obtained through equivariant embeddings of 
flag manifolds. The asymptotics of the Kronecker coefficients are then governed by the properties of the normal bundle of the embedding, whose weights must be contained in an open half space. When this occurs, a very general stability phenomenon 
can be observed, and the limit multiplicities can be computed. In particular the local structure of the Kronecker 
polyhedron can be described (Theorem \ref{general}). 
\item The simplest example is a Segre embedding of a product of projective spaces. This takes care of the very first instance 
of stability, discovered by Murnaghan a very long time ago \cite{mur1,mur2}, which corresponds to  
the simplest stable triple $(1,1,1)$.
This immediately leads to an expression of the limit multiplicities, traditionnally called {\it reduced Kronecker coefficients}. 
Moreover this can readily be generalized: considering the product of a projective space by a Grassmannian we recover  
a  stability property called $k$-{\it stability} by Pak and Panova \cite{pp}, and we are able to compute the limit 
mutiplicities (Proposition \ref{kstab}). Of course this further generalizes to a general product of Grassmannians, which yields a new stability property. 
In fact, what we show is that $(1^{ab}, a^b, b^a)$ is always a stable triple  (Proposition \ref{grass}). 
Moreover we can in principle compute the limit multiplicities. 
\item Turning to products of flag varieties we explain why the convexity property of the normal bundle of an embedding 
defined by a standard tableau is in fact equivalent to the property that the tableau is additive  (Proposition \ref{stab}). 
In fact this is already 
contained in \cite{man-if}, even in its version generalized to multitableaux. What we explain in the present paper is how the 
additive tableaux define minimal regular faces of the Kronecker polyhedron  (Proposition \ref{minimalfaces}). 
Moreover, we describe all the facets 
containing these minimal faces in terms of special tableaux that we call maximal relaxations  (Proposition \ref{relax}). The corresponding inequalities 
are completely explicit  (Propositions \ref{inequalities} and \ref{inequalities2}).
\item Each of these  special facets consists in stable triples, that we therefore obtain in abundance. Moreover, for each
of these triples there is a corresponding notion of reduced Kronecker coefficients. When the stable triple is regular, which
is the generic case, we compute this reduced Kronecker coefficient as a number of integral points in an explicit 
polytope  (Proposition \ref{polytope}). As we mentionned above, this is something we would very much like to be able to do for general 
Kronecker coefficients. 
\item The last section of the paper is devoted to certain Kronecker coefficients for partitions of rectangular shape.  
There are some nice connections with the beautiful theory of $\theta$-groups of Vinberg and Kac  (Proposition \ref{finite}). 
In particular we prove an identity suggested  by Stembridge 
in \cite{stem} by relating it to the affine Dynkin diagram of type $E_6$ (see Proposition \ref{affine}, which also contains 
an interesting identity coming from affine $D_4$). 
\end{enumerate}

\noindent {\it Acknowledgements.} This paper was begun in Berkeley during the semester on  Algorithms and Complexity in Algebraic Geometry organized at the Simons Institute for Computing, and completed in Montr\'eal at the Centre de Recherches Math\'ematiques (Universit\'e de Montr\'eal) and the CIRGET  (UQAM). The author 
warmly thanks these institutions for their generous hospitality. 

\section{From Kronecker to Borel-Weil}

\subsection{Schur-Weyl duality}
For any integer $n$, the irreducible complex representations of the symmetric group $S_n$ are naturally parametrized by 
partitions of $n$ \cite{man-book}. We denote by $[\lambda]$ the representation defined by the partition $\lambda$ of $n$
(we use the notation $\lambda\vdash n$ to express that $\lambda$ is a partition of $n$, in which case we also say that 
$\lambda$ has size $n$). The Kronecker
coefficients can be defined as dimensions of spaces of $S_n$-invariants inside triple tensor products:
$$ g(\lambda,\mu, \nu) = \dim (  [\lambda]\otimes [\mu]\otimes [\nu])^{S_n}.$$
Since the irreducible representations of $S_n$ are defined over the reals, they are self-dual and therefore, 
we can also interprete the Kronecker coefficients as multiplicities in tensor products:
 $$  [\lambda]\otimes [\mu]= \bigoplus_{\nu\vdash n}  g(\lambda,\mu, \nu) [\nu].$$
Schur-Weyl duality allows to switch from representations of symmetric groups to representations of general linear groups. 
Recall that irreducible polynomial representations of $GL(V)$ are parametrized by partitions with at most $d=\dim V$ non 
zero parts. We denote by $S_\lambda V$ the representation defined by the partition $\lambda$. Schur-Weyl duality
can be stated as an isomorphism
$$ V^{\otimes n} = \bigoplus_{\lambda\vdash n} [\lambda]\otimes S_\lambda V $$
between $S_n\times GL(V)$-modules. Here the summation is over all partitions $\lambda$ of $n$ with at most $d$ non 
zero parts (the number of non zero parts will be called the length and denoted $\ell (\lambda)$). 
A straightforward consequence is that 
$$  g(\lambda,\mu, \nu) = \mathrm{mult} ( S_\lambda V \otimes S_\mu W,  S_\nu (V\otimes W)),$$ 
as soon as the dimensions of $V$ and $W$ are large enough; more precisely, as soons as $\dim(V)\ge\ell(\lambda)$
and $\dim(W)\ge\ell(\mu)$. Note that consequently, we get the classical result that 
$$  g(\lambda,\mu, \nu) =0 \qquad \mathrm{if}\; \ell(\nu)>\ell(\lambda)\ell(\mu).$$

\subsection{The Borel-Weil theorem and its consequences}
 The $GL(U)$-representation  $S_\lambda U$ is called a Schur module. It is defined more generally for $\lambda$ an arbitrary non increasing sequence of integers of length equal to the dimension of $U$.  The Borel-Weil theorem asserts that such a 
Schur module can be realized as a space of sections of a linearized line bundle $L_\lambda$ on the complete flag variety 
$Fl(U)$. We can therefore understand the Kronecker coefficients as the mutiplicities in the decomposition of 
$$  S_\nu (V\otimes W)=H^0(Fl (V\otimes W), L_\nu)$$ 
into irreducible $GL(V)\times GL(W)$-modules.

A direct (well-known) consequence of the Borel-Weil theorem is that the direct sum of all the Schur powers of a given 
vector space $U$ has a natural algebra structure. This algebra is finitely generated, and therefore if we let $U=V\otimes W$,
the subalgebra of $GL(V)\times GL(W)$-covariants is also finitely generated (this is a consequence of the fact that the 
unipotent radical of a reductive algebraic group is a Grosshans subgroup). This implies the following result. 
Consider 
$$Kron_{a,b,c}:=\{(\lambda,\mu, \nu),\; \ell(\lambda)\le a,  \ell(\mu)\le b,  \ell(\nu)\le c, \:
g(\lambda,\mu, \nu)\ne 0\}.  $$

\begin{prop}[Semigroup property]\label{sg}
$Kron_{a,b,c}$ is a finitely generated semigroup. 
\end{prop}

(The semigroup property is Conjecture 7.1.4 in \cite{kly}, where finite generation is also a conjecture.)
Moreover, since the covariant algebra is described in terms of saces of sections of line bundles, it has no zero divisor.
This implies  the following monotonicity property: 

\begin{prop}[Monotonicity]\label{mon}
If  $g(\lambda,\mu, \nu)\ne 0$, then for any triple    $(\alpha,\beta, \gamma)$,
 $$g(\alpha+\lambda,\beta+\mu, \gamma+\nu)\ge g(\alpha,\beta, \gamma).$$
\end{prop}

\medskip\noindent {\bf Remark.}
In fact a stronger property is true: there exists a natural map
$$  (  [\lambda]\otimes [\mu]\otimes [\nu])^{S_n}\otimes (  [\alpha]\otimes [\beta]\otimes [\gamma])^{S_p}
\longrightarrow (  [\lambda+\alpha]\otimes [\mu+\beta]\otimes [\nu+\gamma])^{S_{n+p}} $$
which is non zero on  decomposable tensors. It could be interesting to understand this map better.

\subsection{The Kronecker polyhedron}

The semigroup $Kron_{a,b,c}$ lives inside ${\bf Z}^{a+b+c}$, more precisely inside a codimension two sublattice
because of the obvious condition $|\lambda|=|\mu|=|\nu|$ for a Kronecker coefficient $  g(\lambda,\mu, \nu)$ to be non zero. 
Consider the cone generated by  $Kron_{a,b,c}$. The finite generation of the latter implies that this is a rational polyhedral 
cone  $PKron_{a,b,c}$, defined by some finite list of linear inequalities. 

The interpretation in terms of sections of line bundles on flag manifolds allows to understand this polytope as a 
moment polytope, and to use the powerful results that has been obtained in this context using the tools of 
Geometric Invariant Theory, or in an even broader, not necessarily algebraic context, those of Symplectic Geometry. 

On the GIT side, there exist general statements that allow, in principle, to determine moment polytopes and 
in particular the Kronecker polyhedron   $PKron_{a,b,c}$. For example Klyachko in \cite{kly} suggested to use 
the results of Berenstein and Sjamaar, and applied them in small cases. This method has recently been improved in \cite{vw}. This approach has two important limitations. First, it 
describes the moment polytope by a collection of inequalities which is in general far from being minimal: this means
that even if one was able to list these inequalities, many of them would in fact not correspond to any facet of the 
polytope and would be redundant with the other ones. Second, even if far from being satisfactory, the effective 
production of this collection of inequalities would require the computation of certain Schubert constants in some complicated 
homogeneous spaces, which seems combinatorially an extremely hard challenge. More precisely, one would need to
decide whether certain Schubert constants are zero or not, which is certainly more accessible than computing them
but is certainly far beyond our current level of undertanding. 

The first issue has been essentially solved by Ressayre \cite{ressayre}, who devised a way to produce, 
in principle, a
minimal list of inequalities of certain moment polytopes: that is, a list of its facets, or codimension one faces. This 
approach has been remarkably successful for Littlewood-Richardson coefficients. However,
applying them concretely to the Kronecker polyhedron   $PKron_{a,b,c}$ seems completely untractable at the moment. 

At least do we know the Kronecker polyhedron   $PKron_{a,b,c}$ for small values of $a,b,c$. For example (\cite{kly, franz})
$PKron_{3,3,3}$  is defined by the following seven inequallities, and those obtained by permuting the partitions 
$\lambda, \mu, \nu$:
$$\begin{array}{rcl}
\lambda_1+\lambda_2 & \le & \mu_1+\mu_2+ \nu_1+\nu_2, \\
\lambda_1+\lambda_3 & \le & \mu_1+\mu_2+ \nu_1+\nu_3, \\
\lambda_2+\lambda_3 & \le & \mu_1+\mu_2+ \nu_2+\nu_3, \\
\lambda_1+2\lambda_2 & \le & \mu_1+2\mu_2+ \nu_1+2\nu_2, \\
\lambda_2+2\lambda_1 & \le & \mu_1+2\mu_2+ \nu_2+2\nu_1, \\
\lambda_3+2\lambda_2 & \le & \mu_1+2\mu_2+ \nu_3+2\nu_2, \\
\lambda_3+2\lambda_2 & \le & \mu_2+2\mu_1+ \nu_2+2\nu_3.
\end{array}$$

Klyachko made extensive computations showing an overwhelming growth of complexity when the parameters increase. For example 
he claims that $PKron_{2,3,6}$ is defined by $41$ inequalities, $PKron_{2,4,8}$ by  $234$ inequalities and 
 $PKron_{3,3,9}$ by no less than $387$ inequalities!

\subsection{The quasipolynomiality property}

The problem of understanding multiplicities in spaces of sections of line bundles 
is at the core of the study of Hamiltonian actions on symplectic manifolds. 
The {\it Quantization commutes with Reduction} type results have very strong consequences in the context we are interested in.
In the algebraic setting, one starts with a smooth projective complex variety $M$ with an action of a reductive group $G$,
and an ample  $G$-linearized line bundle $L$ on $M$. In this context one considers the virtual $G$-module
$$RR(M,L)=\oplus_{q\ge 0}(-1)^q H^q(M,L),$$
 whose dimension is given by the Grothendieck-Riemann-Roch formula. It is then a very general statement, due to 
 Meinrenken and Sjamaar \cite{ms} (see also \cite{vergne}), that for any dominant weight $\alpha$, the multiplicity of
 the irreducible $G$-module of highest weight $k\alpha$ inside the virtual $G$-module $RR(M,L^k)$ is given 
by a quasipolynomial function of $k$. 

Of course we will apply this result to $G=GL(U)\times GL(V)\times GL(W)$ acting on $M={\bf P}(U\otimes V\otimes W)$, and 
$L$ the hyperplane line bundle on this projective space.  Since ample line bundles on projective spaces have no higher
cohomology, we are in the favorable situation where $RR(M,L^k)$ is just the actual $G$-module $H^0(M,L^k)=S^k(U\otimes V\otimes W)^*$, whose multiplicities are precisely the Kronecker coefficients. Applying the full force of 
Meinrenken and Sjamaar's results, we deduce the following statement:

\begin{theo}  
The stretched Kronecker coefficient  $g(k\lambda,k\mu, k\nu)$ is a piecewise quasi-polynomial function of  
$(k,\lambda,\mu, \nu)$. More precisely,  there is a finite decomposition of the Kronecker polyhedron 
$PKron_{a,b,c}$ into 
closed polyhedral subcones called chambers, and for each chamber $C$ a quasi-polynomial $p_C(k,\lambda,\mu, \nu)$, such that 
$$g(k\lambda,k\mu, k\nu) = p_C(k,\lambda,\mu, \nu)$$
whenever $(\lambda,\mu, \nu)$ belongs to $C$.
\end{theo}

\begin{coro}
For any triple $(\lambda,\mu, \nu)$, the stretched Kronecker coefficient  $g(k\lambda,k\mu, k\nu)$ is a  
quasi-polynomial function of $k\ge 0$. 
\end{coro}

\medskip 
The monotonicity property of Kronecker coefficients easily implies that whenever $(\lambda,\mu,\nu)$ is in the interior of the Kronecker polyhedron, the  stretched Kronecker coefficient  $g(k\lambda,k\mu, k\nu)$ grows as 
fast as possible, in the sense that  if $g(\lambda,\mu, \nu)\ne 0$, then 
$$g(k\lambda,k\mu, k\nu)\simeq  C(\lambda,\mu,\nu) k^{n_{gen}}.$$
Here we denoted by $n_{gen}$ the generic order of growth, which is also the generic dimension 
of the GIT-quotient $M//G$. For $PKron_{a,b,c}$ we get the value 
$$ n_{gen} = abc-a^2-b^2-c^2+2.$$
Moreover the coefficient $C(\lambda,\mu,\nu)$ can be expressed as the  volume of the so-called reduced space  $M_{(\lambda,\mu, \nu)}$ (see \cite{vergne}, page 21).
This volume function is given by the Duistermaat-Heckman measure,
which is piecewise polynomial, and not just quasipolynomial. 

Conversely, if $g(\lambda,\mu, \nu)\ne 0$ 
and $g(k\lambda,k\mu, k\nu)$ grows like $k^n$ for some $n<n_{gen}$, 
then the triple  $(\lambda,\mu, \nu)$ must belong to the boundary of $PKron_{a,b,c}$. The extreme case is when 
$n=0$. As observed in \cite{stem}, this can happen only if  $g(k\lambda,k\mu, k\nu)=1$ for all $k\ge 1$. 
Otherwise said, $(\lambda,\mu, \nu)$ is a {\it weakly stable triple}. 

\smallskip
 The quasi-polynomiality property has attracted the attention of several authors. Most of them realized that general arguments 
based on finite generation of covariant rings implied the eventual quasipolynomiality of the stretched Kronecker coefficients, that is, 
$g(k\lambda,k\mu, k\nu)$ is given by some quasipolynomial for large enough $k$. The much stronger property that the 
 stretched Kronecker coefficients are quasipolynomial right from the beginning seem to be much less accessible using 
algebraic methods only. It was asked as a question in \cite{stem}. It is also discussed in \cite{mul} and some explicit
 computations  appear in \cite{bor1}.

\medskip\noindent
{\bf Remark 1.} One interesting  implication of the quasipolynomiality property is that, knowing the 
Kronecker coefficients asymptotically, in fact we know them completely. Could this be useful to find triples $(\lambda,\mu, \nu)$
in the Kronecker polyhedron, 
such that $g(\lambda,\mu, \nu)=0$? The study of such {\it holes} seems extremely challenging, but is of the greatest importance
for Geometric Complexity Theory \cite{blmw}. 

\medskip\noindent
{\bf Remark 2.} We have noticed that although the  stretched Kronecker coefficient  $g(k\lambda,k\mu, k\nu)$ is a  
quasi-polynomial function of $k\ge 0$, its highest order term is really polynomial. More generally, examples 
show that the period of the coefficients seem to increase when one consider terms of lower degrees. It would be interesting 
to prove an explicit statement of this type.

\section{Kronecker coefficients and Taylor expansions}

\subsection{The general setup}
In \cite{man-if} we suggested a method which easily produces lots of stable triples, and yields much more informations about Kronecker coefficients. The idea is very general: we can study a space  
of sections of any line bundle $L$ on any smooth irreducible variety $X$ by taking the Taylor expansion of these sections along a smooth subvariety $Y$. To be more precise, we can define a filtration of $H^0(X,L)$ by the order of vanishing along $Y$; more formally, we  let  
$$ F_i=  H^0(X,I_Y^i\otimes L) \subset H^0(X,L),$$
where $I_Y$ is the ideal sheaf of $Y$. Let $\iota : Y\hookrightarrow X$ denote the inclusion map. The quotients $I_Y^i/I_Y^{i+1}=\iota_*S^iN^*$, for $N=\iota^*TX/TY$ the normal vector bundle of $Y$ in $X$. Therefore 
there are natural injective maps
$$ F_i/F_{i+1} \hookrightarrow  H^0(Y,\iota^*L\otimes S^iN^* ).$$
In words, this map takes a section of $L$ vanishing to order $i$ on $Y$, to the degree $i$ part of its Taylor expansion in the 
directions normal to $Y$. The injectivity is clear: if the degree $i$ part of the Taylor expansion is zero, then the section 
vanishes on $Y$ at order $i+1$. 

In general, it will be quite difficult to determine the image of these maps. But the situation improves dramatically if 
we suppose that $L$ is a sufficiently large tensor power of some given ample line bundle $M$ on $X$. Indeed, it 
follows from the general properties of ample line bundles that for any fixed integer $i$, the map 
 $$ F_i/F_{i+1} \hookrightarrow  H^0(Y,\iota^*M^k\otimes S^iN^* )$$
must be surjective for $k$ large enough. (This is a straightforward and very classical consequence of 
Serre's vanishing theorems for ample sheaves.)

Let us suppose moreover that the whole setting is preserved by the action of some reductive group $G$. 
Then our filtration  splits as a filtration by $G$-submodules, and a splitting yields an injection of $G$-modules
$$  H^0(X,M^k) \hookrightarrow H^0(Y,\iota^*M^k\otimes S^*N^* ).$$
In fact this statement certainly holds true without any ampleness assumption: it simply asserts that an algebraic section of a line bundle is completely determined by its full Taylor expansion. In case $M$ is ample, we have a very important extra information: 
we know that the right hand side is generated by the left hand side up to any given degree, if 
$k$ is sufficiently large. 

\subsection{Application to Kronecker coefficients}
We want to apply the previous ideas in the following situation. The variety $X$ will be a flag manifold of a tensor product
$V\otimes W$, not necessarily a complete flag; its type (by which we mean the sequence of dimensions of the subspaces 
in the flag) will be allowed to vary, so we will just denote it by $Fl_*(V\otimes W)$. The line 
bundle $M$ will be some $L_{\lambda}$. If we need it to be ample the jumps in the partition $\lambda$ will have 
to be given exactly by the type of the flag manifold, in which case we will say that $\lambda$ is strict (one should add: 
relatively to the flag manifold under consideration). The subvariety $Y$ of $X$ will be a product $Fl_*(U)\times 
Fl_*(V)$ of flag manifolds of certain types, and we will need to construct the embedding $\iota$. Of course our 
reductive group will be $G=GL(V)\times GL(W)$. 

We have $Y=G/P$ for some parabolic subgroup $P$ of $G$. It is well known that the category of $G$-equivariant 
vector bundles on $Y$ is equivalent to the category of finite dimensional $P$-modules. Such modules can be quite 
complicated since $P$ is not reductive. Let us consider a Levi decomposition $P=LU_P$, where $U_P$ denotes the
unipotent radical of $P$ and $L$ is a Levi factor, in particular a reductive subgroup of $P$. We can consider $L$-modules
as $P$-modules with trivial action of $U_P$, and conversely. In general, any $P$-module has a canonical filtration 
by $P$-modules, such that the associated graded $P$-module has a trivial action of $U_P$. We can then consider it
as an $L$-module and decompose it as a sum of irreducible modules, as in the usual theory of representations of 
reductive groups. Equivalentely, any $G$-homogeneous vector bundle $E$ as an associated graded homogeneous
bundle $gr(E)$, which is a direct sum of irreducible ones. Cohomologically speaking, and even at the level of global
sections, $E$ and $gr(E)$ can be quite different. In fact $E$ can in principle be reconstructed from $gr(E)$ through 
a series of extensions than can mix up the cohomology groups in a complicated fashion. Note however that $H^0(Y,E)$ 
is always a $G$-submodule of  $H^0(Y,gr(E))$.

\medskip\noindent {\bf Example.}
On a complete flag variety, the irreducible bundles are exactly the line bundles. On a variety of incomplete 
flags  $Fl_*(U)$ this is no longer true, but we can obtain them as follows.  
Denote by $p: Fl(U)\rightarrow Fl_*(U)$ the natural projection map
(which takes a complete flag and simply forgets some of the subspaces in it). Consider a line bundle 
$L_\lambda$ on  $Fl(U)$. Then the pushforward  
$E_\lambda=p_*L_\lambda$ is an irreducible homogeneous vector bundle on   $Fl_*(U)$ and they are all obtained 
in this way. (In particular, if $\lambda$ is strict then $E_\lambda$ is exactly the line bundle $L_\lambda$ on $Fl_*(U)$, 
for which we use the same notation.) Note that the Borel-Weil theorem extends to all the irreducible bundles:
$$ H^0( Fl_*(U), E_\lambda)  = H^0( Fl(U), L_\lambda)  = S_\lambda U. $$

\medskip\noindent
{\bf Definition}. The embedding $\iota$ is {\it stabilizing} if the normal bundle has the following convexity 
property: the highest  weights of the positive degree part of 
the symmetric algebra $S^*(gr(N^*))$, considered as an $L$-module, are contained in some open half space of the 
weight lattice of $L$.
 
\medskip
The line bundle $\iota^*L_\lambda$ will be of the form $L_{a(\lambda)}\otimes L_{b(\lambda)}$ for some weights 
$a(\lambda)$ and  $b(\lambda)$ depending linearly on $\lambda$. More generally, given an irreducible homogeneous 
vector bundle $E_\alpha$ on $Fl_*(V\otimes W)$, we will need to understand its pull-back  by $\iota$. The resulting 
homogeneous bundle will scarcely be completely reducible. The associated graded bundle is of the form
$$  gr(\iota^*E_\alpha)=\bigoplus_{(\rho,\sigma)\in T(\alpha)}E_\rho\otimes E_\sigma$$ 
for some multiset $T(\alpha)$.

 Our main result  is the following. 

\begin{theo}\label{general}
Suppose that the embedding $\iota$ is stabilizing. Let $\lambda$ be any strict partition. Then 
\begin{enumerate}
\item The Kronecker coefficient $g(k\lambda, ka(\lambda), kb(\lambda))=1$ for any $k\ge 0$. 
\item For any triple $(\alpha, \beta, \gamma)$ the Kronecker coefficient  $g(\alpha+k\lambda, \beta+ka(\lambda), 
\gamma+kb(\lambda))$ is a non decreasing, bounded function of $k$, hence eventually constant.
Otherwise said, the triple$(\lambda, a(\lambda), b(\lambda))$ is stable. 
\item If moreover $\lambda$ is ample, then the limit Kronecker coefficient  is given by the multiplicity of the weight 
$(\beta, \gamma)$ inside $gr(\iota^*E_\alpha)\otimes S^*(gr(N^*))$. 
\item In particular, if $\alpha$ is strict, this is 
the multiplicty of $(\beta-a(\alpha), \gamma-b(\alpha))$ inside $S^*(gr(N^*))$.
\end{enumerate}
\end{theo}

\noindent {\it Proof.} The Kronecker coefficients  $g(k\lambda, ka(\lambda), kb(\lambda))$ are positive since the restriction map
$$ S_{k\lambda}(V\otimes W)= H^0( Fl_*(U), L_\lambda^k)\longrightarrow  H^0( Fl_*(V)\times Fl_*(W), \iota^*L_\lambda^k)
= S_{ka(\lambda)}V\otimes S_{kb(\lambda)}W$$
is surjective. Indeed this map is non zero since the line bundle $ L_\lambda$  is generated by global sections, 
which means that it has a section that 
does not vanish at any given point. The surjectivity then follows from Schur's lemma since the right hand side is irreducible. 
Moreover, we have seen that Taylor expansions induces a $G$-embedding
$$ S_{k\lambda}(V\otimes W) \hookrightarrow  H^0( Fl_*(V)\times Fl_*(W), \iota^*L_\lambda^k\otimes S^*(N^*)).$$
As we also noticed, replacing $N^*$ by  $gr(N^*)$ can only result in making the right hand side larger, so we even 
have an inclusion
$$ S_{k\lambda}(V\otimes W) \hookrightarrow  H^0( Fl_*(V)\times Fl_*(W), \iota^*L_\lambda^k\otimes S^*(gr(N^*))).$$
By the Borel-Weil theorem, the right hand side is now a sum of irreducible $GL(V)\times GL(W)$-modules whose  highest 
weights are of the form $(ka(\lambda),kb(\lambda))$ plus a highest weight of  $ S^*(gr(N^*))$. But since these weights are
supposed to belong to a strictly convex cone,  none can give a positive contribution of highest weight 
 $(ka(\lambda),kb(\lambda))$ except the weight zero in degree zero. This proves that 
 $g(k\lambda, ka(\lambda), kb(\lambda))=1$ for any $k\ge 0$.

Now consider  an arbitrary triple $(\alpha, \beta, \gamma)$. By the Borel-Weil theorem again, we have 
$$ S_{\alpha+k\lambda}(V\otimes W) = H^0(Fl_*(U),E_\alpha\otimes L_\lambda^k).$$
This immediately implies that $g(\alpha+k\lambda, \beta+ka(\lambda), 
\gamma+kb(\lambda))$ is a non decreasing function of $k$, by considering the map 
 $$ H^0(Fl_*(U),E_\alpha\otimes L_\lambda^k)\otimes H^0(Fl_*(U),L_\lambda) 
\longrightarrow H^0(Fl_*(U),E_\alpha\otimes L_\lambda^{k+1})$$
 and its restriction through $\iota$. Moreover,
the same approach using Taylor expansions can then be used without much change, in particular we get an injection 
$$ S_{\alpha+k\lambda}(V\otimes W) \hookrightarrow  H^0( Fl_*(V)\times Fl_*(W), gr(\iota^*E_\alpha)\otimes
\iota^*L_\lambda^k\otimes S^*(gr(N^*))).$$
The same argument as before then implies that $g(\alpha+k\lambda, \beta+ka(\lambda))$ is bounded by the 
multiplicity of the representation of highest weight $(\beta,\gamma)$ inside 
$ gr(\iota^*E_\alpha)\otimes S^*(gr(N^*))$, which is finite by the convexity hypothesis. 

What remains to check is that for $k$ large enough, we have in fact equality. But being finite, the multiplicity 
$(\beta,\gamma)$ inside $ gr(\iota^*E_\alpha)\otimes S^*(gr(N^*))$ only comes from some finite part 
of the symmetric algebra of the conormal bundle. Then it follows formally from the properties of ample line 
bundles that if $k$ is large enough, all the maps involved in the Taylor expansions are surjective in boundeed degrees, 
and it makes no difference to replace all the homogeneous bundles involved by the associated graded bundles. 
This concludes the proof. $\Box $

\begin{coro}\label{stable}
The set of triples $(\lambda, a(\lambda), b(\lambda))$ is a face of the Kronecker polyhedron. Moreover the local structure 
of the polyhedron around this face is given by the cone generated by the weights of  $S^*(gr(N^*))$. 
\end{coro}

\begin{coro}\label{genmon}
Limit Kronecker coefficients have the monotonicity property.  
\end{coro}

In particular, if they are non zero, then they remain nonzero after stretching. Kirillov \cite{kir} and Klyachko \cite{kly} conjectured in certain special cases that the converse property, that is saturation, 
 should hold. We could ask more generally if our  limit 
(or reduced) Kronecker coefficients should have the saturation property, but we actually have only very 
limited evidence for that.

\medskip
\noindent
{\bf Remarks.}
\begin{enumerate}
\item We have not tried to bound the minimal value of $k$ starting from which the Kronecker coefficient  
 $g(\alpha+k\lambda, \beta+ka(\lambda), \gamma+kb(\lambda))$ does stabilize, but it would in principle be possible to give an effective version of the previous theorem. Indeed an irreducible homogeneous bundle which has non zero sections 
has no higher cohomology by Bott's theorem, and this would be sufficient to ensure for example that replacing 
the conormal bundle by its graded associated bundle makes no difference at the level of global sections. Nevertheless the 
resulting statements would probably be rather heavy, and presumably  not even close to be sharp. 
\item We have expressed the limit multiplicities in a rather compact form, but we will see that these expressions are in fact very complicated. Moreover they usually involve other Kronecker coefficients, as well as Littlewood-Richardson coefficients and 
other interesting variants. It seems we are playing with Russian dolls, but with growing complexity when we open a doll to 
find a smaller, more mysterious one... 
\end{enumerate}

\subsection{Tangent and normal bundles}

In order to apply the previous theorem, we need to be able to understand the normal bundle of our embedding 
$\iota$ and the associated graded bundle. The starting point for that is of course to understand the tangent 
bundle to a flag variety; this is classical and goes as follows. 
 
Suppose we consider a variety $Fl_*(V)$ of flags $(0=V_0\subset V_1\subset \cdots \subset V_{r-1}\subset V_r=V)$
where $V_i$ has dimension $d_i$. Each of these spaces defines a homogeneous vector bundle on  $Fl_*(V)$, but not 
irreducible in general. The irreducible homogeneous bundles are obtained by considering the quotient bundles 
$Q_i=V_i/V_{i-1}$, for $1\le i\le r$. Each of these bundles is irreducible, as well as each of their Schur powers. 
More generally, each irreducible vector bundle on $Fl_*(V)$ is of the form
$$  E_\lambda = S_{\alpha_1}Q_1\otimes \cdots \otimes  S_{\alpha_r}Q_r$$
for some non increasing sequences $(\alpha_1,\ldots ,\alpha_r)$ of relative integers. Each homogeneous vector bundle 
$F$ can then be constructed from such irreducible bundles by means of suitable extensions; the set of irreducible 
bundles involved does not depends on the process and their sum is $gr(F)$. For example we have non trivial 
exact sequences $0\rightarrow V_{i-1}\rightarrow V_{i}\rightarrow Q_{i}\rightarrow 0$, and by induction we deduce that 
$$gr(V_i)=Q_1\oplus\cdots \oplus Q_i.$$ 
The tangent bundle to     $Fl_*(V)$ at the flag  $V_\bullet=(0=V_0\subset V_1\subset \cdots \subset V_{r-1}\subset V_r=V)$
can be naturally identified with the quotient of $End(V)$ by the subspace of endomorphisms preserving   $V_\bullet$. 
Taking the orthogonal with respect to the Killing form we get a natural identification of the cotangent bundle:
$$ T^*_{Fl_*(V)} \simeq \{ X\in End(V), \; X(V_i)\subset V_{i-1}, 1\le i\le r\}.$$
From this we can easily deduce the following statement:

\begin{lemm}\label{tangent}
The associated bundle of the tangent bundle of  $Fl_*(V)$ is
$$ gr( T_{Fl_*(V)} ) = \oplus_{1\le i<j\le r}Hom(Q_i,Q_j).$$
\end{lemm}

In order to describe the  normal bundle of an embedding $\iota : Fl_*(V)\times Fl_*(W)\hookrightarrow FL_*(V\otimes W)$, 
we will have to take the quotients of two such bundles and it is clear that the result will be quite complicated in general. 
Rather than writing down  general formulas, we will compute these quotients in certain specific situations; mainly, when there
are only few terms (typically for Grassmannians, whose tangent bundles are irreducible) or for complete flag varieties 
(whose irreducible bundles are line bundles, hence much easier to handle). 

\section{First examples}

We now have all the necessary ingredients in hand in order to apply Theorem \ref{general}.  
What remains to be done is to construct suitable embeddings between flag varieties 
and determine whether they are stabilizing. We begin with a few simple examples.

\subsection{Murnaghan's stability}
As an appetizer we begin with a well-known instance of stability, first discovered by Murnaghan \cite{mur1,mur2}.
This is the statement that for any triple of partitions $(\alpha, \beta ,\gamma)$, the Kronecker coefficient
$g(\alpha+(k), \beta+(k), \gamma+(k))$ is eventually constant for large $k$ (a sharp bound on $k$ has been given in \cite{bor3}). 
The limit value is called a reduced Kronecker coefficient and is denoted  $\bar{g}(\alpha', \beta' ,\gamma')$, where 
$\alpha'$ is deduced from $\alpha$ by suppressing the first part. 

To interpret this in our setting, consider the Segre embedding 
$$\iota : {\bf P}(V)\times {\bf P}(W)\hookrightarrow {\bf P}(V\otimes W).$$ 
If $L$ denotes the tautological line bundle on projective space, then with some abuse of notations, 
$\iota^*L=L\otimes L$. From Lemma \ref{tangent} we easily derive that 
$$  N=  Hom(L,V/L)\otimes Hom(L,W/L). $$
In particular this embedding is clearly  stabilizing since any component of $S^*N$ will have negative degree on both
copies of $L$. This implies Murnaghan's stability right away.

Moreover, if $Q$ is the quotient bundle on ${\bf P}(V\otimes W)$, then the associated bundle of $\iota^*Q$ is 
$$  gr(\iota^*Q) = L\otimes W/L \oplus V/L\otimes L\oplus V/L\otimes W/L.$$
Applying Theorem \ref{general}, we deduce that $\bar{g}(\alpha', \beta' ,\gamma')$ is  the 
multiplicity of $S_{\beta' }A\otimes S_{\gamma' }B$ inside 
$$S_{\alpha' }(A\otimes B\oplus A\oplus B)\otimes S^*(A\otimes B).$$
This is equivalent to Lemma 2.1 in \cite{bor3}. The monotonicity property that is a special case of
Corollary  \ref{genmon}
has been observed in \cite{chris}. 

\subsection{$k$-stability}

A generalisation of Murnaghan's stability has been considered recently by Vallejo \cite{vall-diag}
and Pak and Panova \cite{pp} 
This $k$-stability, as coined by the latter authors,  corresponds to the generalized Segre embedding 
$$\iota : {\bf P}(V)\times Gr(k,W)\hookrightarrow Gr(k,V\otimes W),$$ 
sending a pair $(L,M)$ made of a one dimensional subspace $L\subset V$ and a $k$ dimensional subspace $M\subset W$, 
to the $k$ dimensional subspace $L\otimes M$ of  $V\otimes W$.
If $L$ denotes the tautological line bundle on the Grassmann 
variety, then $\iota^*L=L^k\otimes L$. Extending the computation we made for  projective spaces we get that   
$$  gr N= Hom(L,V/L)\otimes End_0(M) \oplus
Hom(L,V/L)\otimes Hom(M,W/M). $$
For essentially the same reasons as in the previous case, this embedding is manifestly stabilizing. 
This implies Theorem 1.1 in \cite{pp} without further ado. Theorem 10.2 in \cite{vall-diag} gives an
effective version.

Moreover we have access to the limit multiplicities, called the $k$-reduced Kronecker coefficients. 
For this we need to consider a vector bundle $E_\lambda$ on 
 $Gr(k,V\otimes W)$ and pull it back by $\iota$. Since $E_\lambda$ is a Schur power of a the quotient bundle $Q$ 
and $$gr(\iota^*Q)=L\otimes W/M \oplus V/L\otimes W\oplus V/L\otimes W/M,$$
we get the following result.

\begin{prop}\label{kstab}
The $k$-reduced Kronecker coefficient $\bar{g}_k(\lambda,\mu,\nu)$ is equal to the multiplicity of the product 
 $A_-^{\mu_-}\otimes S_{\mu_+}A_+\otimes S_{\nu_-}B_-\otimes S_{\nu_+}B_+$ inside
$$S_\lambda( A_-\otimes B_+\oplus A_+\otimes B_-\oplus A_+\otimes B_+)\otimes Sym(A_-^*\otimes A_+\otimes End_0(B_-))
\otimes Sym(A_-^*\otimes A_+\otimes B_-^*\otimes B_+). $$
\end{prop}

For Geometric Complexity Theory the most relevant Kronecker coefficients  are those indexed by two partitions with equal
rectangular shapes \cite{blmw}. This corresponds to taking $\lambda$ and $\nu$ equal to the empty partition. To get a contribution 
from the previous formula we must avoid all the terms contributing positively to $B_+$, which kills the second symmetric 
algebra. The first will contribute through $S_\mu A_+\otimes S_\mu(End_0(B_-))$ (because of the Cauchy formula), and 
extracting the term with $\nu_-=0$ means that we take the $GL(B_-)$-invariants inside $S_\mu(End_0(B_-))$. Since 
$End_0(B_-)$ is just $sl_k$, we get: 

\begin{coro}
The $k$-reduced Kronecker coefficient $\bar{g}_k(0,\mu,0)$ is equal to the dimension of the $GL_k$-invariant subspace of 
$S_\mu (sl_k)$.
\end{coro}

An effective version of this result was derived in \cite{man-kr} using a completely different method. 
\smallskip

It seems interesting to notice that Schur powers of $sl(V)=End_0(V)$ are often involved in the stable Kronecker coefficients, 
along with the Kronecker coefficients themselves and the Littlewood-Richardson coefficients. These multiplicities are also of interest for themselves; more generally, many interesting phenomena appear when one considers Schur powers of the adjoint 
representation of a simple complex Lie algebra, see eg \cite{kostant}. Our methods can be applied to the study of the 
asymptotics of these coefficients; we hope to come back to this question in a future paper.

\subsection{Grassmannian stability}
We can obviously generalize the Segre embedding to any product of Grassmann varieties.
For any positive integers $a,b$ consider the natural embedding  
$$\iota : Gr(a,V)\times Gr(b,W)\hookrightarrow Gr(ab,V\otimes W).$$ 
If $L$ denotes the tautological line bundle on the Grassmann variety, then $\iota^*L=L^b\otimes L^a$. Extending the computation we made for  projective spaces we get that   
$$  gr N= End_0(A)\otimes Hom(B,W/B) \oplus Hom(A,V/A)\otimes End_0(B) \oplus
Hom(A,V/A)\otimes Hom(B,W/B). $$
This embedding is again stabilizing and we get right away the following generalization of $k$-stability. 

\begin{prop}\label{grass} Let $a,b$ be positive integers. 
For any triple $(\lambda,  \mu , \nu)$, the Kronecker coefficient 
$$g(\lambda+(t)^{ab},  \mu+ (at)^b, \nu+(bt)^a)$$
 is a non decreasing, eventually constant, function of $t$. 
\end{prop}

We could call the limits $(a,b)$-reduced Kronecker coefficients, and provide an expression generalizing Theorem 3. 
We leave this to the interested reader. Note that in this statement  $(\lambda,  \mu , \nu)$ are not just partitions 
but integer sequences for which the arguments of the Kronecker coefficient are partitions, at least for large enought $t$. 
For example $\lambda=(\lambda_+,\lambda_-)$ where $\lambda_+$ is a non decreasing sequence of length $ab$ (with
possibly negative entries), and  $\lambda_-$ is a partition (of arbitrary length). 

Note that when $\lambda_-$ is the empty partition, drastic simplifications occur. This follows from the usual formula 
$S_{\alpha+(t)^d}V=S_\alpha V\otimes (\det V)^t$ when $d$ is the dimension of $V$. Since $\det (A\otimes B)=
(\det A)^{\dim B}\otimes (\det B)^{\dim A}$, we deduce that 
$$g(\lambda+(t)^{ab},  \mu+ (at)^b, \nu+(bt)^a)=g(\lambda,  \mu, \nu).$$
This is Theorem 3.1 in \cite{vallejo}.

\section{Stability and standard tableaux}

\subsection{Classification of embeddings}

In this section we focus on the equivariant embeddings of complete flag varieties 
$$ \iota : Fl(A)\times Fl(B)\hookrightarrow Fl(A\otimes B),$$
the combinatorics being more transparent in that case. Denote by $a$ and $b$ the respective dimensions of $A$ and $B$.
The following two statements appear in \cite{man-if}:

\begin{prop}
The equivariant embeddings $ \iota : Fl(A)\times Fl(B)\hookrightarrow Fl(A\otimes B)$ are classified by 
standard tableaux of rectangular shape $a\times b$.
\end{prop}

Let $T$ be such a standard   tableau of rectangular shape $a\times b$. The corresponding embedding $\iota_T$ is defined as follows. 
Let $A_\bullet$ and $B_\bullet$ be two complete flags in $A$ and $B$ respectively. Choose 
adapted basis $u_1,\ldots ,u_a$ of $A$ and 
 $v_1,\ldots ,v_b$ of $B$, that is, such that $A_i=\langle  u_1,\ldots ,u_i\rangle$ and 
 $B_j=\langle  v_1,\ldots ,v_j\rangle$. Then define the complete flag $C_\bullet$ in $A\otimes B$ by 
 $$C_k=\langle u_i\otimes v_j, \; T(i,j)\le k\rangle.$$
Here $T(i,j)$ denotes the entry of $T$ in the box $(i,j)$. Clearly this flag does not depend on the adapted basis but only on the flags, and we can let $\iota_T(A_\bullet, B_\bullet)=C_\bullet$.

\begin{prop}\label{stab}
The  embedding $ \iota_T : Fl(A)\times Fl(B)\hookrightarrow Fl(A\otimes B)$ is stabilizing if and only if the  
standard tableau $T$ is additive.
\end{prop}

Here we use the terminology of Vallejo, the additivity condition being used in \cite{vallejo2}.  In our paper \cite{man-if} no special 
terminology was introduced for this additivity property, which means the following: there exist increasing sequences 
$x_1<\cdots <x_a$ and $y_1<\cdots <y_b$ such that 
$$ T(i,j)<T(k,l) \quad\Longleftrightarrow \quad  x_i+y_j<x_k+y_l.$$

\noindent {\bf Remark.} 
 There is a huge number of embeddings $\iota_T$: recall that by the hook length formula, the number of standard tableaux of shape $a\times b$ is $ ST(a,b)= (ab)!/h(a,b)$, where $h(a,b)=(a+b-1)!!/(a-1)!!(b-1)!!$. This grows at least exponentially with $a$ and $b$. Among these, the proportion of additive tableaux probably tends to zero when $a$ and $b$ grow, but their number should still 
grow exponentially.  Note that each of this additive tableau corresponds to a certain chamber in the complement 
of the arrangement of hyperplanes $H_{ijkl}$ defined by the equalities   $x_i+y_j=x_k+y_l$  in ${\bf R}_+^{a+b-2}$ 
(we may suppose that $x_1=y_1=0$). This looks very much like the hyperplance arrangement associated to a root system. 

\medskip
The partitions $a_T(\lambda)$ and  $b_T(\lambda)$  such that $\iota_T^*L_\lambda = L_{a_T(\lambda)}\otimes  L_{b_T(\lambda)}$
are easily described; one just needs to read  the entries in each line or column of $T$ and sum the corresponding parts of $\lambda$: 
$$ a_T(\lambda)_i=\sum_{j=1}^b \lambda_{T(i,j)}, \qquad  b_T(\lambda)_j=\sum_{i=1}^a \lambda_{T(i,j)}.$$

\subsection{$(T,\lambda)$-reduced Kronecker coefficients}

Applying our general statements of Theorem \ref{general}, we get:

\begin{prop}\label{stabletriples}
Let $T$ be any additive tableau.
For  any partition $\lambda$, the triple $(\lambda, a_T(\lambda), b_T(\lambda))$
is stable. 
\end{prop}

Although this statement does not appear explicitely in \cite{man-if}, it is discussed p.735, in the paragraph 
just before the Example. It was recently rediscovered by Vallejo \cite{vallejo2}, inspired by the work of 
Stembridge \cite{stem}, and following a completely different approach. 
(As a matter of fact  there is a slight difference between our definition of additivity and that of Vallejo.
When the tableau $T$ is additive, as we intend it, then for any partition $\lambda$ of length at most $ab$, the matrix $A$ with entries $a_{ij}= \lambda_{T(i,j)}$ is additive in the sense of Vallejo, and all such matrices 
are obtained that way. As we just stated it, this implies that the triple $(\lambda, a_T(\lambda), b_T(\lambda) )$ 
is stable. If moreover $\lambda$ is strict, then we will get more information, in particular we will be able 
to compute the stable Kronecker coefficients. )

Note that the additivity property is introduced in section 3.1.2 of \cite{man-if} and explained to be equivalent
to the convexity property of the normal bundle that implies stability. At that time we were mainly interested in plethysm 
and we treated the case of Kronecker coefficients rather quickly, giving details only for a sample of the results 
that were amplified in the strongly similar case of plethysm. Also we treated directly the multiKronecker 
coefficients
$$g(\mu_1, \ldots ,\mu_r) =\dim ([\mu_1]\otimes \cdots \otimes [\mu_r])^{S_n}, $$
for which the method applies with essentially no difference. 
\smallskip

The fact that additivity is equivalent to stability is easy to understand. Recall, as a special case of Lemma \ref{tangent},  that the 
tangent   bundle of a complete flag manifold $Fl(U)$ has associated graded bundle
$$gr(T_{Fl(U)})=\oplus_{i<j}Hom(Q_i,Q_j),$$ 
where the quotient bundles are now line bundles. Applying this to $U=V\otimes W$ and pulling-back by $\iota_T$ we will get a formula in terms of the quotient line bundles on $Fl(V)$ and $Fl(W)$, that we will denote by $E_i$ and $F_j$. The formula reads
$$gr(\iota_T^*T_{Fl(V\otimes W)})=\oplus_{ T(i,j)<T(k,l)}Hom(E_i\otimes F_j,E_k\otimes F_l).$$
Let us denote by $e_1,\ldots ,e_a$ and $f_1,\ldots ,f_b$ natural basis of the weight lattices of $GL(V)$ and $GL(W)$, respectively. The formula shows that   the weights of the restricted tangent 
bundle $\iota_T^*(T_{Fl(V\otimes W)})$ are the $e_k+f_l-e_i-f_j$ for $ T(i,j)<T(k,l)$. 
We claim that the weights of the normal bundle are exactly the same. Indeed, since the tableau $T$ is 
increasing on rows and columns, the weights  $e_k-e_i$ and $f_l-f_j$ appear for $k>i$ and $l>j$ with 
multiplicity $b$ and $a$ respectively, in particular greater than one. Since they are the weights of the tangent bundles of $Fl(V)$ and $Fl(W)$, 
we see that when we go from the restricted tangent bundle to the normal bundle the list of weights 
will not change, only certain multiplicities will decrease, but remaining positive. In particular the generated cone 
will not be affected. 

Finally, recall that the additivity property asks for the existence of sequences  $x_1 <\cdots <x_a$ and $y_1 <\cdots <y_b$
such that  $ T(i,j)<T(k,l)$ if and only if $x_i+y_j<x_k+y_l$, otherwise written as $x_k-x_i+y_l-y_j>0$. 
This precisely means that the corresponding linear form 
is positive on each of the weights  $e_k-e_i+f_l-f_j$ with $ T(i,j)<T(k,l)$. 

\medskip\noindent {\bf Remark.} Of course it is not 
necessary, in order to determine the cone generated by the weights of $gr(N^*)$, to compile the full list of $M=(a-1)(b-1)(ab+a+b)/2$ weights. The $ab-1$ of them obtained by reading the successive entries of the tableau $T$ from 
$1$ to $ab$ will suffice, since all the other weights will obviously be sums of these.  

\medskip
Note that from Proposition \ref{stabletriples}, we can deduce right away the following special case, which is a generalization of Proposition \ref{grass}.

\begin{coro}
For  any partition $\mu$ of size $m$, the triple $(1^m, \mu, \mu^*)$ is stable.
\end{coro}

\noindent {\it Proof. } Consider a rectangle $a\times b$ in which the diagram of $\mu$ can be inscribed. Consider 
$\mu$ and the conjugate partition $\mu^*$ as integer sequences of lengths $a$ and $b$ respectively, by adding 
zeroes if necessary. Consider the increasing sequences $\alpha_1, \ldots , \alpha_a$ and  $\beta_1, \ldots , \beta_b$
defined by   $\alpha_i=i-\mu_i-1$ and $\beta_j=j-\mu_j^*$. Then $\alpha_i+\beta_j=-h_{ij}$, the opposite of the 
hook length of $\mu$ for the box $(i,j)$. In particular   $\alpha_i+\beta_j$ is negative exactly on the support 
of $\mu$. Let $T$ be the corresponding additive tableau, and let $\lambda=1^m$. Then $a_T(\lambda)=\mu$
and $b_T(\lambda)=\mu^*$. Hence   $(1^m, \mu, \mu^*)$ is stable. $\Box$

\medskip
Let us denote the value of  $ g(\alpha+k\lambda, \beta+ka_T(\lambda), \gamma+kb_T(\lambda))$, for $k$ very large, by 
$g_{T,\lambda}(\alpha, \beta, \gamma)$, and call it a $(T,\lambda)$-{\it reduced Kronecker coefficient}. If $\lambda$ is strictly decreasing of length $ab$, or $ab-1$, so that the corresponding line bundle on the flag manifold is very ample, we have seen that this $(T,\lambda)$-reduced Kronecker coefficient can be computed as the multiplicity of the weight $(\beta-a_T(\alpha), 
\gamma-b_T(\alpha))$ inside the symmetric algebra $S^*(gr(N^*))$. The weights of $gr(N^*)$ can readily be read off the tableau $T$. Let us denote them by $(u_i,v_i)$, for $1\le i\le M=(a-1)(b-1)(ab+a+b)/2$. 
Denote by $P_{T,\lambda}(\mu,\nu)$ the polytope 
defined as the intersection of the quadrant $t_1,\ldots , t_M\ge 0$ in ${\bf R}^M$ with the affine 
linear space defined by the condition that 
$$ \sum_{i=1}^Nt_i(u_i,v_i)=(\mu,\nu).$$

\begin{prop}\label{polytope}
The $(T,\lambda)$-reduced Kronecker coefficient $g_{T,\lambda}(\alpha, \beta, \gamma)$  is equal to the number 
of integral points in the polytope  $P_{T,\lambda}(\beta-a_T(\alpha), \gamma-b_T(\alpha))$. 
\end{prop}

Of course we would then be tempted to stretch the triple $(\alpha, \beta, \gamma)$. We would then obtain the 
stretched $(T,\lambda)$-reduced Kronecker coefficient $g_{T,\lambda}(k\alpha, k\beta, k\gamma)$ as given by 
the Ehrhart quasi-polynomial of the rational polytope $P_{T,\lambda}(\beta-a_T(\alpha), \gamma-b_T(\alpha))$. 
This suggests interesting behaviours for multistretched  Kronecker coefficients, but we will not pursue on this route.

Let us simply notice an obvious consequence for  $(T,\lambda)$-reduced Kronecker coefficients, the following translation invariance property:
$$  g_{T,\lambda}(\alpha, \beta, \gamma)=g_{T,\lambda}(\alpha+\delta, \beta+a_T(\delta), \gamma+b_T(\delta))$$
for any partition $\delta$. 

\subsection{Faces of the Kronecker polytope}

As we observed, the convexity property of the embeddings defined by additive tableaux has very interesting consequences for the Kronecker polytope. Let us first reformulate Corollary \ref{stable}.

\begin{prop}\label{minimalfaces} 
Each additive tableau $T$ defines a regular face $f_T$ of the Kronecker polyhedron $PKron_{a,b,ab}$, of minimal dimension.
\end{prop}
 
Regular means that the face meets the interior of the Weyl chamber, which is the set of strictly decreasing partitions. 
Ressayre proved in \cite{ressayre}, in a more general setting, that the maximal codimension of a regular face is the rank
of the group which in our case is $GL(V)\times GL(W)$. This exactly matches with the codimension of $f_T$. 
  
Around this minimal face $f_T$, we know that the  local structure of  $PKron_{a,b,ab}$
  is described by the convex polyhedron generated by the weights of $gr(N^*)$. 

\smallskip
A face of our polytope will be defined by a linear function which is non negative on all these vectors, and vanishes on 
a subset of them that generate a hyperplane. Such a linear function will be defined by sequences 
$x_1 \le \cdots \le x_a$ and $y_1 \le\cdots \le y_b$ such that  $x_i+y_j\le x_k+y_l$ when $ T(i,j)<T(k,l)$. 
The different values of  $x_i+y_j$ define a partition of the rectangle $a\times b$ into disjoint regions, and this partition 
is a {\it relaxation} of $T$, in the sense that each region is numbered by consecutive values of $T$. 
The hyperplane condition can be interpreted as the fact that the vectors  $e_k-e_i+f_l-f_j$, for $(i,j)$ and $(k,l)$
belonging to the same region, generate a hyperplane in the weight space. This is also a maximality condition: we 
cannot relax any further while keeping the compatibility condition with $T$. We deduce:

\begin{prop}\label{relax}
The facets $F_R$ of the Kronecker polyhedral cone $Kron_{a,b,ab}$ containing the minimal face $f_T$ are in
bijective correspondence with the maximal relaxations $R$ of the tableau $T$. 
\end{prop}

\noindent {\bf Example.} Recall that for $a=b=3$ there are $42$ standard tableaux fitting in a square of size three, among which $36$ are additive. The total number of maximal relaxations of these additive tableaux is $17$. Up to diagonal symmetry they are as follows, where each square is filled by the entry $x_i+y_j$ in box $(i,j)$ for some sequences  $(x_1=0,x_2,x_3)$ and $(y_1=0,y_2,y_3)$.

$$\begin{array}{ccccc}
\begin{array}{ccc}
0&0&0 \\ 0&0&0 \\ 1&1&1
\end{array} &
\begin{array}{ccc}
0&0&0 \\ 1&1&1 \\ 1&1&1
\end{array} &
\begin{array}{ccc}
0&0&1 \\ 1&1&2 \\ 2&2&3
\end{array} &
\begin{array}{ccc}
0&1&1 \\ 1&2&2 \\ 2&3&3
\end{array} &
\begin{array}{ccc}
0&1&2 \\ 2&3&4 \\ 3&4&5
\end{array} \\
 &&&& \\
\begin{array}{ccc}
0&0&1 \\ 0&0&1 \\ 1&1&2
\end{array} &
\begin{array}{ccc}
0&0&1\\ 1&1&2 \\ 2&2&3
\end{array} &
\begin{array}{ccc}
0&1&1 \\ 1&2&2 \\ 1&2&2
\end{array} &
\begin{array}{ccc}
0&1&2 \\ 1&2&3 \\ 2&3&4
\end{array} &
\begin{array}{ccc}
0&1&2 \\ 1&2&3 \\ 3&4&5
\end{array}
\end{array}$$
Consider for example the relaxation $R$ encoded in the tableau 
$$\begin{array}{ccc}
0  &1&2 \\ 1&2&3 \\ 3&4&5
\end{array}$$
It splits the square into six regions, three of size one and three of size two. There are therefore eight compatible standard tableaux $T$, which are all additive. This gives eight minimal faces $f_T$ incident  to the facet $F_R$.

\begin{prop}\label{inequalities}
The defining inequalities of the facet $F_R$ associated to a maximal relaxation $R$ defined by sequences 
$(x_i,y_j) $ are of the form 
$$ \sum_{i=1}^ax_i\beta_i+  \sum_{j=1}^by_j\gamma_j\ge  
\sum_{i=1}^a   \sum_{j=1}^b(x_i+y_j)\alpha_{T(i,j)},$$
where $T$ is any standard tableau compatible with $R$. 
\end{prop}

We mean that $R$ is a relaxation of $T$. It is clear then that the inequality does not depend on $T$, since taking another 
compatible tableau $T'$ amounts to switching some entries only inside the regions defined by $R$, on each of which the 
sum   $x_i+y_j$ is constant.

\subsection{From rectangles to arbitrary tableaux}
If we want to restrict to the Kronecker cone  $Kron_{a,b,c}$ for some $c<ab$, we just need to intersect 
 $Kron_{a,b,ab}$ with a linear space $L_c$ of codimension $ab-c$, which meets each of our minimal faces $f_T$. 
But this raises two issues. First, two distinct minimal faces can have the same intersection with $L_c$. Second, 
it is not clear whether the intersection of a facet $F_R$ of  $Kron_{a,b,ab}$ with $L_c$ will still be a facet of 
 $Kron_{a,b,c}$.   

 The first issue is easy to address: two tableaux $T$ and $T'$ will give the same minimal face of  $Kron_{a,b,c}$
if and only if they coincide up to $c$, that is, their entries  smaller of equal to $c$ appear in the same boxes. We thus
get  minimal faces of  $Kron_{a,b,c}$ parametrized by standard tableaux $S$ of size $c$ inside the rectangle $a\times b$, 
which are additive in the same sense as before.   

To address the second issue we can modify our embeddings $\iota_T$ accordingly. A standard tableau $S$ of size $c$ inside the rectangle $a\times b$ defines an embedding 
$$ \iota_S : FL(V)\times Fl(W)\hookrightarrow FL_c(V\otimes W),$$
where we denote by $FL_c(U)$ the partial flag manifold of $U$ parametrizing flags of the form $(0=U_0\subset U_1
\subset \cdots \subset U_c\subset U)$, where $U_i$ has dimension $i$. We will require that $S$ fits exactly in the rectangle, not in any smaller one.   Denote the quotient bundles on $Fl(U)$ by $Q_1, \ldots , 
Q_c, Q_{c+1}$, they are all line bundles except the last one. We have $\iota_S^*Q_k=E_i\otimes F_j$ whenever
$S(i,j)=k\le c$, while  $\iota_S^*Q_{c+1}$ is not irreducible, but has associated graded bundle  
 $$ gr(\iota_S^*Q_{c+1})=\oplus_{(i,j)\notin S}E_i\otimes F_j.$$
The pull-back of the tangent bundle is then given by the same formula as before, 
 $$gr(\iota_S^*T_{Fl(V\otimes W)})=\oplus_{ S(i,j)<S(k,l)}Hom(E_i\otimes F_j,E_k\otimes F_l),$$
except that for this to be correct, we need to consider $S$ as a rectangular tableau of size $a\times b$, in which 
the boxes $(i,j)$ that do not belong to $S$ are all numbered by the same arbitrarily large number, say $S(i,j)=\infty$.  
Exactly as before we deduce that the weights of the normal bundle are the   $e_k-e_i+f_l-f_j$ with $ S(i,j)<S(k,l)$,
and the convexity condition translates into the same additive property, that we can summarize by saying that 
$S$ must be a piece of a rectangular additive standard tableau. The discussion above then goes through exactly as in the 
rectangular case, except that we don't need to care about the boxes of the rectangle that are not supported by $S$. 
We get the following slight extension of our previous results: 

\begin{prop}\label{inequalities2}
Let $S$ by a standard tableau of height $a$, width $b$, size $c$. Suppose that $S$ is additive. Then the  set of 
stable triples of the form $(\lambda, a_S(\lambda), b_S(\lambda))$ defines a minimal regular face $f_S$ of the 
Kronecker polytope $PKron_{a,b,c}$. Moreover the facets of the polytope containing this minimal face $f_S$ are in bijection 
with the minimal relaxations $R$ of $S$. If $R$ is defined by non decreasing sequences $(x_i,y_j)$, 
the equation of the facet $F_R$ is
$$ \sum_{i=1}^ax_i\beta_i+  \sum_{j=1}^by_j\gamma_j\ge  
\sum_{(i,j)\in S} (x_i+y_j)\alpha_{S(i,j)}.$$
\end{prop}

\medskip\noindent {\bf Remark.}
More variants could be explored. Maps $Fl_*(V)\times Fl_*(W)\hookrightarrow Fl_*(V\otimes W)$ 
where arbitrary types appear at the source are easily constructed in terms of tableaux, and their stability could be analyze.
This will be more complicated in general since the quotient bundles will have ranks bigger than one. Moreover it will only 
give access to stable triples on the boundary of the Weyl chamber.

Of course we could aso readily extend the discussion
to an arbitrary number of vector spaces, describe embeddings of arbitrary products of flag manifolds in terms of 
multidimensional tableaux, observe that the convexity condition on the weights of the normal bundle is again an additivity
condition, and deduce stability properties for multiKronecker coefficients. This was partly done in \cite{man-if}.

\section{Rectangles: stability and beyond}

There exist stable triples which do not come from additive tableaux, and it would be nice to understand them. 
Stembridge in \cite{stem} observed that $(22,22,22)$ is stable, and it is certainly not additive. (Nevertheless it is highly degenerate, in the sense that it belongs to a very small face of the dominant Weyl chamber. In particular this example leaves open the question of the existence of non additive stable triples in the interior of the Weyl 
chamber.) In this section
we make a connection with Cayley's hyperdeterminant and the Dynkin diagram $D_4$, and we explain another 
observation by Stembridge in terms of affine $E_6$. 

\subsection{Finite cases}

The rectangular Kronecker coefficients, ie those involving partitions of rectangular shape, are of special interest
because of their direct relation with invariant theory. For three factors, 
$$ g(p^a,q^b,r^c)=\dim (S^k(A\otimes B\otimes C))^{SL(A)\times SL(B)\times SL(C)} $$
when $k=pa=qb=rc$ and $a,b,c$ are the dimensions of $A, B, C$. 
 Of course this connection also holds for a larger number of factors. 

There seems to exist only few results on these quotients. The cases for which there are only finitely many orbits of 
$GL(A)\times GL(B)\times GL(C)$ inside $A\otimes B\otimes C$  have been completely classified in connection with 
Dynkin diagrams(see eg \cite{man-prehom} and references therein). 
One can deduce all the possible dimensions $a,b,c$ for which there exists a dense orbit, through
the combinatorial process called castling transforms. In this situation there exists one invariant
for $SL(A)\times SL(B)\times SL(C)$ at most, depending on the
codimension of the complement of the dense open orbit, which can be one (in which case its equation is an invariant)
or greater than one (in which case there is no invariant).  The former case gives a weakly stable triple. 

In this classification through Dynkim diagrams, the triple tensor products correspond to triple nodes, so there are
only few cases, coming from diagrams of type D or E. The first interesting case is $D_4$, corresponding to $a=b=c=2$.
The invariant is the famous hyperdeterminant first discovered by Cayley, which has degree four. This implies that 
$g(n^2,n^2,n^2)=1$ when $n$ is even and $g(n^2,n^2,n^2)=0$ when $n$ is odd. In particular, $(22,22,22)$ is a
weakly stable triple, and even a stable triple, as shown by Stembridge \cite{stem}. 

The next two diagrams, $D_5$ and $D_6$, give  $(a,b,c)=(2,2,3)$ and $(2,2,4)$
respectively. The corresponding invariants have degree $6$ and $4$. They correspond to the weakly stable triples 
$(33,33,222)$ and $(22,22,1111)$ (which we already met among stable triples). 
For $D_n$, $n\ge 7$, we get  $(a,b,c)=(2,2,n-2)$ but there is no non trivial invariant anymore. 
Finally the triple nodes of $E_6, E_7, E_8$ yield the triples $(a,b,c)=(2,3,3), (2,3,4), (2,3,5)$. There is no invariant for 
the latter case, but an invariant of degree $12$ in the two previous ones, yielding the weakly stable triples $(66,444,444)$
and $(66,444,3333)$. Let us summarize our discussion: 

\begin{prop}\label{finite}
The triple nodes of the Dynkin diagrams of types $D_4$, $D_5$, $E_6$, $E_7$ yield the non additive weakly stable triples 
$(22,22,22)$, $(33,33,222)$, $(66,444,444)$, $(66,444,3333)$.
\end{prop}

\subsection{Affine cases}

This discussion can be upgraded from Dynkin to affine Dynkin diagrams. Indeed it is a theorem of Kac \cite{kac}
that when we consider a
representation associated to a node of an affine Dynkin diagram, the invariant algebra is free, or otherwise said, is a 
polynomial algebra. (If the chosen node is the one that has been attached to the usual Dynkin diagram, the associated representation is just the adjoint one, so  the theorem generalizes the well-known result that the invariant algebra 
of the adjoint representation is free.) 

There will be four cases related to Kronecker coefficients, corresponding to the affine Dynkin diagrams with a unique multiple node. 
Let us introduce the following notation:
$$\begin{array}{rcl}
g_{\hat{D}_4}(n) & = & g(n^2,n^2,n^2,n^2), \\
g_{\hat{E}_6}(n) & = & g(n^3,n^3,n^3), \\
g_{\hat{E}_7}(n) & = & g((2n)^2,n^4,n^4), \\
g_{\hat{E}_8}(n) & = & g((3n)^2,(2n)^3,n^6).
\end{array}$$
Applying Kac's results we immediately identify the generating series of these rectangular  Kronecker coefficients.

\begin{prop}\label{affine}
The generating series of the rectangular  Kronecker coefficients 
$g_{\hat{D}_4}(n)$,
$g_{\hat{E}_6}(n)$, 
$g_{\hat{E}_7}(n)$ and 
$g_{\hat{E}_8}(n)$ are the following:
$$\sum_{n\ge 0}g_{\hat{D}_4}(n)q^n  =  \frac{1}{(1-q)(1-q^2)^2(1-q^3)}, $$
$$\sum_{n\ge 0}g_{\hat{E}_6}(n)q^n =  \frac{1}{(1-q^2)(1-q^3)(1-q^4)}, $$
$$\sum_{n\ge 0}g_{\hat{E}_7}(n)q^n  =   \frac{1}{(1-q^2)(1-q^3)}, $$
$$\sum_{n\ge 0}g_{\hat{E}_8}(n)q^n  =   \frac{1}{1-q}.
$$
\end{prop}

The last of these identities simply expresses the fact that $(33,222,1^6)$ is a weakly stable triple, as we already know. 
The previous one can be rewritten as 
$$ g((2n)^2,n^4,n^4)=\frac{n+\pi_6(n)}{6}$$
where $\pi_6$ is the 6-periodic function with first 6 values $(6,-1,4,3,2,1)$. The identity for $g_{\hat{E}_6}(n)$ 
has been suggested by Stembridge (\cite{stem}, Appendix). It can also be rewritten as a quasipolynomial identity: 
$$  g(n^3,n^3,n^3) = \frac{(n+1)(n+2)}{48}+\frac{n+1}{16}\pi_2(n)+\frac{1}{48}\pi_{12}(n),$$
where $\pi_{12}$ is 12-periodic with period $(37,-12,9,16,21,-48,25,0,21,4,9,0)$  and
 $\pi_2$ is 2-periodic with period $(3,1) $. Finally the quadruple Kronecker coefficient 
$$    g(n^2,n^2,n^2,n^2) = \frac{n^3+12n^2+29n+18}{72}+\frac{n+1}{72}\pi_2(n)+\frac{1}{72}\pi_6(n),$$
where $\pi_6$ is 6-periodic with period $(35,-8,27,8,19,0)$  and
 $\pi_2$ is 2-periodic with period $(19,10)$.

\medskip\noindent {\bf Remark.} Let us mention that there should exist lots of non additive stable 
triples. For example,  \cite{man-kr} implies that if $\lambda$ is a partition of size $2n$, then 
$(n^2,n^2,\lambda)$ is weakly stable when $\lambda$ is even and of length at most four, as well as 
$(n^4,(2n)^2,2\lambda)$  when $\lambda$ has  length at most three and $\lambda_1\le 
\lambda_2+\lambda_3$.  It would be interesting  to prove that these weakly stable triples are in 
fact stable, and to get more examples, or more general procedures to construct (weakly) stable triples. 

Theorem 6.1 in \cite{stem} gives a criterion for stability that covers 
the additive triples, but not only those. It would be interesting to understand these non additive 
triples more explicitely, find a geometric interpretation, compute the stable Kronecker coefficients, 
and decide to which extent they could help to describe the Kronecker polyhedra.


\begin{thebibliography}{aa}

\bibitem{bor1}
Briand E.,  Orellana R.,  Rosas M., {\it 
Quasipolynomial formulas for the Kronecker coefficients indexed by two two-row shapes}
 21st International Conference on Formal Power Series and Algebraic Combinatorics, 241-252, 
Discrete Math. Theor. Comput. Sci., Nancy, 2009. 

\bibitem{bor2}
Briand E.,  Orellana R.,  Rosas M., {\it 
Reduced Kronecker coefficients and counter-examples to Mulmuley's strong saturation conjecture}, 
with an appendix by K. Mulmuley,
Comput. Complexity {\bf 18} (2009),  577-600. 

\bibitem{bor3}
Briand E.,  Orellana R.,  Rosas M., {\it 
The stability of the Kronecker product of Schur functions},
J. Algebra {\bf 331} (2011), 11-27. 

\bibitem{blmw}
Buergisser P., Landsberg J.M., Manivel L., Weyman J., {\it 
An overview of mathematical issues arising in the Geometric Complexity Theory approach to VP vs VNP},
SIAM Journal in Computing {\bf 40} (2011), 1179-1209.

\bibitem{chris}
Christandl M., Harrow A.W., Mitchison G.,
{\it Nonzero Kronecker coefficients and what they tell us about spectra}, 
Comm. Math. Phys. {\bf 270} (2007), 575-585.

\bibitem{franz}
Franz M.,  {\it 
Moment polytopes of projective $G$-varieties and tensor products of symmetric group representations},
J. Lie Theory {\bf 12} (2002),  539-549.

\bibitem{fulton}
Fulton W., {\it 
Eigenvalues, invariant factors, highest weights, and Schubert calculus}, 
Bull. Amer. Math. Soc. {\bf 37} (2000),  209-249. 

\bibitem{kac}
Kac V.G., {\it 
Some remarks on nilpotent orbits}, 
J. Algebra {\bf 64} (1980),  190-213. 

\bibitem{kir}
Kirillov A.N., {\it 
An invitation to the generalized saturation conjecture}, 
Publ. Res. Inst. Math. Sci. {\bf 40} (2004), 1147-1239. 

\bibitem{kly}
Klyachko A., {\it Quantum Marginal problem and representations of the symmetric group}, arXiv:quant-ph:0409113.
 
\bibitem{kostant}
Kostant B., {\it Clifford algebra analogue of the Hopf-Koszul-Samelson theorem, 
and the $g$-module structure of $\Lambda g$}, Adv. Math. {\bf 125} (1997), 275-350.

\bibitem{man-if}
Manivel L., {\it 
Applications de Gauss et pl\'ethysme},
Ann. Inst. Fourier {\bf 47} (1997), 715-773. 

\bibitem{man-book}
Manivel L.,
Symmetric Functions, Schubert Polynomials and Degeneracy Loci, SMF/AMS {\bf 6}, 2001.

\bibitem{man-kr}
Manivel L., {\it 
A note on certain Kronecker coefficients}, Proc. A.M.S. {\bf 138} (2010), 1-7. 

\bibitem{man-kr}
Manivel L., {\it 
On rectangular Kronecker coefficients}, 
J. Algebraic Combin. {\bf 33} (2011),  153-162. 

\bibitem{man-prehom}
Manivel L., {\it 
Prehomogeneous spaces and projective geometry}, 
to appear in the Rendiconti del Seminario Matematico, Universita e Politecnico di Torino.

\bibitem{ms}
Meinrenken E., Sjamaar R., {\it 
Singular reduction and quantization}, 
Topology {\bf 38} (1999), 699-762. 

\bibitem{mul}
Mulmuley K., {\it 
Geometric Complexity Theory VI: the flip via saturated and positive integer programming in representation theory and algebraic geometry}, arXiv:0704.0229.

\bibitem{mur1}
Murnaghan F.D., {\it  
The analysis of the Kronecker product of irreducible representations of the symmetric group},
Amer. J. Math. {\bf 60} (1938), 761-784. 

\bibitem{mur2}
Murnaghan F.D., {\it  
On the analysis of the Kronecker product of irreducible representations of $S_n$}, 
Proc. Nat. Acad. Sci. U.S.A. {\bf 41}, (1955). 515-518. 

\bibitem{pp}
Pak I.,  Panova G., {\it  
Bounds on the Kronecker coefficients}, arXiv:1406.2988.

\bibitem{rassart}
Rassart E., {\it 
A polynomiality property for Littlewood-Richardson coefficients},  
J. Combin. Theory Ser. A {\bf 107} (2004), 161-179. 

\bibitem{ressayre}
Ressayre N., {\it  Geometric invariant theory and the generalized eigenvalue problem}, 
Invent. Math. {\bf 180} (2010),  389-441.

\bibitem{stem}
Stembridge J., {\it Generalized stability of Kronecker coefficients},
preprint, August 2014. With an Appendix available on http://www.math.lsa.umich.edu/~jrs/papers

\bibitem{vallejo}
Vallejo E., {\it A stability property for coefficients in Kronecker products of complex $S_n$ characters}, 
Electron. J. Combin. {\bf 16} (2009), no. 1, Note 22, 8 pp.

\bibitem{vall-diag}
Vallejo E., {\it A diagrammatic approach to Kronecker squares}, 
Journal of Combinatorial Theory {\bf127} (2014), 243-285.

\bibitem{vallejo2}
Vallejo E., {\it Stability of Kronecker coefficients via discrete tomography},
arXiv:1408.6219 

\bibitem{vergne}
Vergne M., {\it 
Quantification g\' eom\' etrique et r\' eduction symplectique},  
S\' eminaire Bourbaki, Vol. 2000/2001,
Ast\' erisque {\bf 282} (2002), Exp. 888, viii, 249-278. 

\bibitem{vw}
Vergne M., Walter M.,  {\it Moment cones of representations}, arXiv:1410.8144. 
\end{thebibliography}
\end{document}